\documentclass[11pt]{article}

\usepackage{amscd,amsmath, amssymb}

\newcommand{\la}{\lambda}

\newcommand{\f}{\varphi}

\newcommand{\CC}{\mathbb{C}}

\newcommand{\RR}{\mathbb{R}}
\newcommand{\ZZ}{\mathbb{Z}}

\numberwithin{equation}{section}

\def\eqref#1{(\ref{#1})}
\newcommand{\goth}{\mathfrak}

\newcommand{\arrow}{{\:\longrightarrow\:}}
\newcommand{\Z}{{\Bbb Z}}
\newcommand{\C}{{\Bbb C}}
\newcommand{\R}{{\Bbb R}}

\newcommand{\6}{\partial}
\def\1{\sqrt{-1}\:}

\newcommand{\calo}{{\cal O}}

\renewcommand{\tilde}{\widetilde}
\renewcommand{\bar}{\overline}
\renewcommand{\phi}{\varphi}
\renewcommand{\epsilon}{\varepsilon}
\renewcommand{\geq}{\geqslant}

\newcounter{Mycounter}[section]
\newcounter{lemma}[section]
\newcounter{claim}[section]
\renewcommand{\theclaim}{{Claim \thesection.\arabic{claim}}}
\newcommand{\claim}{%
     \setcounter{claim}{\value{Mycounter}}
     \refstepcounter{claim}
     \stepcounter{Mycounter}
     {\noindent \bf \theclaim:\ }}

\newcounter{sublemma}[section]
\newcounter{corollary}[section]
\newcounter{theorem}[section]
\renewcommand{\thetheorem}{{Theorem \thesection.\arabic{theorem}}}
\newcommand{\theorem}{%
     \setcounter{theorem}{\value{Mycounter}}
     \refstepcounter{theorem}
     \stepcounter{Mycounter}
     {\noindent \bf \thetheorem:\ }}

\newcounter{conjecture}[section]
\newcounter{proposition}[section]
\renewcommand{\theproposition}
       {{Proposition \thesection.\arabic{proposition}}}
\newcommand{\proposition}{%
     \setcounter{proposition}{\value{Mycounter}}
     \refstepcounter{proposition}
     \stepcounter{Mycounter}
     {\noindent \bf \theproposition:\ }}

\newcounter{definition}[section]
\renewcommand{\thedefinition}
       {{Definition~\thesection.\arabic{definition}}}
\newcommand{\definition}{%
     \setcounter{definition}{\value{Mycounter}}
     \refstepcounter{definition}
     \stepcounter{Mycounter}
     {\noindent \bf \thedefinition:\ }}

\newcounter{example}[section]
\newcounter{remark}[section]
\renewcommand{\theremark}{{Remark \thesection.\arabic{remark}}}
\newcommand{\remark}{%
     \setcounter{remark}{\value{Mycounter}}
     \refstepcounter{remark}
     \stepcounter{Mycounter}
     {\noindent \bf \theremark:\ }}

\newcounter{problem}[section]
\newcounter{question}[section]
\makeatletter

\newcommand{\ps@verbit}{%
  \renewcommand{\@oddhead}{%
          \scriptsize
          {LCK manifolds with potential}
          \hfil\tiny {L. Ornea and M. Verbitsky, 14 July 2004 }}
  \renewcommand{\@evenhead}{\@oddhead}
  \renewcommand{\@oddfoot}{\hfil\thepage\hfil}
  \renewcommand{\@evenfoot}{\@oddfoot}}

\@addtoreset{equation}{section}
\@addtoreset{footnote}{section}
\makeatother

\def\blacksquare{\hbox{\vrule width 5pt height 5pt depth 0pt}}
\def\endproof{\hfill\blacksquare}

\begin{document}
\begin{center}
{\LARGE\bf
Locally conformal K\"ahler\\[3mm] manifolds with potential\\[3mm]
}

Liviu Ornea and Misha Verbitsky\footnote{Misha Verbitsky is an EPSRC advanced 
fellow supported by CRDF grant RM1-2354-MO02 and EPSRC grant
GR/R77773/01.

{ {\bf Keywords:} locally 
conformal K\"ahler manifold, K\"ahler potential, 
Lee field, Vaisman manifold, Hopf manifold, Stein variety, 
Poincar\'e-Dulac theorem}

\scriptsize
{\bf 2000 Mathematics Subject 
Classification:} { 53C55, 32G05.}}

\end{center}

{\small 
\hspace{0.15\linewidth}
\begin{minipage}[t]{0.7\linewidth}
{\bf Abstract} \\
A locally conformally K\"ahler (LCK) manifold $M$ is one
which is covered by a K\"ahler manifold $\tilde M$
with the deck transform group acting conformally on
$\tilde M$. If $M$ admits a holomorphic flow, 
acting on $\tilde M$ conformally, it is called 
a Vaisman manifold. Neither the class of LCK
manifolds nor that of Vaisman manifolds is stable
under small deformations. We define a new class
of LCK-manifolds, called LCK manifolds with 
potential, which is closed under small deformations.
All Vaisman manifolds are LCK with potential.
We show that an LCK-manifold with potential
admits a covering which can be compactified
to a Stein variety by adding one point. 
This is used to show that any LCK manifold $M$ with
potential, $\dim M \geq 3$, can be embedded to a Hopf manifold,
thus improving on  similar results for Vaisman manifolds (\cite{ov2}).
\end{minipage}
}

\tableofcontents


\section{Introduction}

\subsection{Motivation}

The study of algebraic geometry of compact non-K\"ahler manifolds
is stymied by the lack of metric structures and unability to apply
most differential geometric and analytic arguments. However, in
many examples of non-K\"ahler manifold, the metric structures
arise naturally. The most common example of these is given
by the locally conformally K\"ahler (LCK) geometry. An LCK manifold
is one admitting a K\"ahler covering $\tilde M$ with the
deck transform acting on $\tilde M$ by conformal transforms.
The first examples of non-K\"ahler manifolds were Hopf manifolds,
that is, the quotients of $\C^n \backslash 0$ by a cyclic
group $\langle A \rangle$ generated by an invertible linear
operator $A:\; \C^n \arrow \C^n$, with all eigenvalues $<1$.
Such a quotient is obviously locally conformally K\"ahler
(see e.g. \cite{go}). For a long time it was conjectured
that all non-K\"ahler complex surfaces are LCK
(this conjecture was disproven in \cite{belgun}). Still, the ubiquity
of LCK manifolds is hard to overestimate.

LCK manifolds have been widely studied in the last 
30 years. They share some properties with the K\"ahler manifolds
(they are stable to blow-up, satisfy an embedding theorem of Kodaira
type, cf. \cite{ov2} 
and the last section here), but they also contrast K\"ahler
geometry. A striking feature is that a deformation of an LCK
manifold is not LCK. Even the smaller (and better understood)
subclass of Vaisman manifolds is not stable at deformations (both
results were proved in \cite{belgun} by finding appropriate
counterexamples). On the other hand, it was proven in \cite{ve} that
the K\"ahler metric of the universal covering  of a Vaisman manifold admits a 
K\"ahler potential. This property is stable to small deformations and hence a
Vaisman manifold deforms to a LCK one (but not necessarily Vaisman).

 This suggests the need for a new subclass of 
LCK manifolds, called LCK manifolds with potential,
stable under deformations and containing Vaisman manifolds. 
It is the aim of this paper to   study the 
LCK manifolds with potential.
We obtain a Kodaira-type embedding theorem for LCK-manifolds
with potential (\ref{_koda_embe_Theorem_}), embedding each into
an appropriate Hopf manifold, not necessarily of Vaisman type.

\subsection{Basic definitions}

Let $(M,I,g)$ be a connected complex Hermitian 
manifold of complex dimension at least $2$. 
Denote by $\omega$ its fundamental Hermitian 
two-form, with the convention $\omega(X,Y)=g(X, IY)$. 

\hfill

\definition 
$(M,I,g)$ is called {\bf locally conformal K\"ahler (LCK)} if there exists a 
\emph{closed} one-form  $\theta$ (called the Lee form) such that:
$$d\omega=\theta\wedge \omega.$$

\hfill

Equivalently, any cover $\tilde M$ of $M$ on which the pull-back
$\tilde\theta$ of $\theta$ is exact (in particular, the universal
cover), carries a K\"ahler form $\Omega =e^{-f}\tilde\omega$,
where $\tilde\theta=df$, and such that
$\pi_1(M)$ acts on $\tilde M$ by holomorphic homotheties. 
Conversely, a manifold admitting such a K\"ahler covering
is necessarily locally conformally K\"ahler.

For many equivalent definitions and examples, the reader is referred
to \cite{drag} and \cite{ov2}. 

The most important subclass of LCK manifolds is
defined by the parallelism of the Lee form with 
respect to the Levi-Civita connection of $g$
(\cite{_Vaisman:Dedicata_}, \cite{drag}).

\hfill

\definition 
An LCK manifold $(M,I,g)$ is called {\bf a Vaisman manifold} if
$\nabla\theta=0$, where $\nabla$ is the Levi-Civita 
connection of $g$.

\hfill

According to \cite{_Kamishima_Ornea_}, the compact Vaisman manifolds
can be defined as an LCK manifolds admitting a conformal
holomorphic flow which does not induce an isometry 
on $\tilde M$.

\hfill

Vaisman geometry is intimately related to Sasakian one (see
\cite{ov1}). On the other hand, the K\"ahler form $\Omega$ on the
universal covering of a Vaisman manifold admits a K\"ahler potential
(cf. \cite{ve} and the next section).


\section{LCK manifolds with potential}


\subsection{LCK manifolds and monodromy}

The main object of this paper is described in the following

\hfill

\definition \label{_LCK_w_p_Definition_}
An {\bf LCK manifold with  potential} is a manifold
which admits a K\"ahler covering $(\tilde M, \Omega)$, and a positive
smooth function $\f:\; \tilde M \rightarrow \RR^{>0}$ 
(the {\bf LCK potential}) satisfying the following conditions:
\begin{enumerate}
\item  It is proper, \emph{i.e.} its level sets are compact.
\item The monodromy map $\tau$ acts on it by multiplication with a constant   
$\tau (\f)=const \cdot \f$.
\item It is a  K\"ahler potential, \emph{i.e.} 
$-\1 \6\bar\6 \f = \Omega$.
\end{enumerate}

\hfill

\remark
Let $\Gamma$ be the monodromy (deck transform) group 
of $\tilde M \rightarrow M$, and
$\chi:\; \Gamma \arrow \R^{>0}$ a homomorphism mapping 
$\gamma\in \Gamma$ to the number $\frac{\gamma(\phi)}{\phi}$.
Consider the covering $\tilde M':= \tilde M/\ker \chi$.
Clearly, $\phi$ defines a K\"ahler potential $\phi'$ on 
$\tilde M'$, and $\phi'$ satisfies assumptions
(1)-(3) of \ref{_LCK_w_p_Definition_}. 
Therefore, we may always choose the covering 
$\tilde M$ in such a way that $\chi:\;\Gamma \arrow \R^{>0}$
is injective.
Further on, we shall always assume 
that $\chi$ is injective whenever an LCK manifold
with potential is considered.

\hfill

\remark
Any closed complex submanifold of an LCK-manifold with potential is also
an LCK-manifold with potential.

\hfill

Note that, by \cite[Pr. 4.4]{ve}, all compact Vaisman manifolds do
have a potential. Indeed, on $\tilde M$, we have $\theta=d\f$, where 
$\phi= \Omega(\theta, I(\theta))$. Using 
the parallelism, one shows that  $\Omega=e^{-\f}\tilde\omega$ is a
K\"ahler form with potential $\f$.

The two non-compact, non-K\"ahler, LCK manifolds constructed
recently by J. Renaud in \cite{renaud} do also have a potential.

In both cases the monodromy is isomorphic to $\ZZ$, hence,
by  \ref{_LCK_phi_proper_Proposition_} below, condition
(1) in the definition is fulfilled. 

\hfill

\claim\label{_chi_discrete_for_LCK-pot_Claim_}
Let $M$ be an LCK-manifold with potential,
$\tilde M$ the corresponding cover with the monodromy
$\Gamma$, and $\chi:\; \Gamma \arrow \R^{>0}$
the injective character constructed above.
Then the image of $\chi$ is discrete in $\R^{>0}$,
in other words, $\Gamma$ is isomorphic to $\Z$.

\hfill

\noindent{\bf Proof:} Let $\gamma$ be a non-trivial element. Then
$\phi$ defines a continuous and proper map 
\[ \tilde M/\langle \gamma\rangle
   \stackrel{\phi_1}\arrow \R^{>0}/\langle \phi(\gamma)\rangle\cong S^1.
\]
Since $\R^{>0}/\langle \phi(\gamma)\rangle$ is compact,
and $\phi_1$ is proper, $\tilde M/\langle \gamma\rangle$ 
is also compact. On the other hand, $\Gamma/\langle \gamma\rangle$
acts freely on $\tilde M/\langle \gamma\rangle$. Therefore,
$\Gamma/\langle \gamma\rangle$ is finite. \endproof

\hfill

The converse is also true, hence:

\hfill

\proposition\label{_LCK_phi_proper_Proposition_}
Let $M$ be an LCK-manifold, $\tilde M \arrow M$ its K\"ahler
covering, $\chi:\; \Gamma \arrow \R^{>0}$
the character defined above, and
 $\phi$ a K\"ahler potential on $\tilde M$ 
which satisfies $\gamma(\phi) = \chi(\gamma) \phi$.
Assume that $\chi$ is injective. Then $\phi$ is proper
if and only if the image of $\chi$ is discrete in
$\R^{>0}$.

\hfill

\noindent{\bf Proof:} The ``only if'' part is proven in
\ref{_chi_discrete_for_LCK-pot_Claim_}.
Assume, conversely, that the image of $\chi$ is discrete.
Then $\Gamma$ is isomorphic to $\Z$,  and $\phi$ defines
a map 
\[ \tilde M/\langle \gamma\rangle
   \stackrel{\phi_1}\arrow \R^{>0}/\langle \phi(\gamma)\rangle,
\]
where $\gamma$ is the generator of $\Gamma$. Since 
$\tilde M/\langle \gamma\rangle=M$, and $M$ is compact,
this map is proper.
Its fibers are naturally identified with the level sets of $\phi$,
and therefore they are compact, and $\phi$ is proper.
\endproof

\subsection{Deformations of LCK manifolds with potential}

With an argument similar to the one in  \cite[Th. 4.5]{ov2}, 
noting that small deformations of positive potentials remain 
positive, we derive

\hfill

\theorem
The class of compact LCK manifolds with potential is stable under
small deformations.

\hfill

\noindent{\bf Proof:} Let $(M,I')$ be a small deformation of $(M,I)$.
Then $\phi$ is a proper function on $(\tilde M, I')$ satisfying
$\gamma(\phi) = \chi(\gamma)\phi$. It is strictly plurisubharmonic 
because a small deformation of a strictly plurisubharmonic
function is again strictly plurisubharmonic (note that the fundamental domain of 
the monodromy action is compact). Therefore,
$(\tilde M, I')$ is K\"ahler, and 
$\phi$ is an LCK-potential on $(\tilde M, I')$. 
\endproof

\hfill

Hence, in particular, any compact Vaisman manifold deforms to a LCK
one. This explains \emph{a posteriori} why the construction in
\cite{go} worked for deforming the Vaisman structure of a Hopf
surface of K\"ahler rank 1 to  a non-Vaisman 
LCK structure on Hopf surface of K\"ahler rank 0 which, moreover, by
\cite{belgun}, does not admit any Vaisman metric.

On the other hand, not all LCK manifolds admit an LCK potential. For
example, one may consider the blow up of a point on a compact
Vaisman manifold. Moreover, the LCK structure of the Inoue surfaces
do not admit potential, since they can be deformed to the non-LCK
type Inoue surface $S^+_{n;p,q,r,u}$ with $u\in \CC\setminus \RR$
(cf. \cite{belgun}).


\section{Holomorphic contractions on Stein varieties}


\subsection{Filling a K\"ahler covering of an LCK-manifold with potential}

\theorem\label{_filling_Theorem_}
Let  $M$ be a LCK manifold with a potential, $\dim M \geq 3$,
and $\tilde M$ the corresponding covering. Then $\tilde M$
is an open subset of a Stein variety $\tilde M_c$, with at most
one singular point. Moreover, $\tilde M_c \backslash \tilde M$
is just one point.

\hfill

\noindent{\bf Proof:}
Consider
the subset 
\[ \tilde M(a) = \{x\in \tilde M\ \  | \ \ \f(x) >a\}.
\]
 By the properties of the potential function, this set  
is holomorphically concave. By Rossi-Andreotti-Siu theorem 
(cf. \cite[Th. 3, p. 245]{rossi} and \cite[Pr. 3.2]{andreotti_siu}) 
$\tilde M(a)$ can be compactified: it is an open subset of 
a Stein variety $\tilde M_c$, with 
(at most) isolated singularities\footnote{It is here that the restriction on the dimension is essential.}. The ring
of holomorphic functions on the "filled in" (see also \cite{ma}) part
of $\tilde M(a)$ is identified with the ring
of CR-functions on the level set $\f^{-1}(a)$.
Therefore, one could extend the embedding 
$\tilde M(a)\rightarrow \tilde M_c$ to 
$\tilde M\rightarrow \tilde M_c$, filling all $\tilde M$
to a Stein variety with at most isolated singularities.

We show that $\tilde M_c$ is obtained from $\tilde M$ by adding only one 
point. 

The deck transform (monodromy) 
group $\Gamma\cong \Z$ acts on $\tilde M_c$
by non-trivial conformal automorphisms. Let
$\gamma$ be its generator, which is a contraction\footnote{It has to be a contraction or an expansion, we make the first choice.}
(that is, $\gamma$ maps $\Omega$ to $\lambda \cdot\Omega$ with $\lambda<1$).
Consider a holomorphic function $f$ on $\tilde M_c$.
The variety $\tilde M_c$ is covered by a sequence
of compact subsets $B_r= \tilde M_c\backslash \tilde M(r)$.
Since $\gamma$ maps $B_r$ to $B_{\lambda r}$, the sequence
$\gamma^n f$ in uniformly bounded on compact sets.
Therefore, it has a subsequence converging
to a holomorphic function, in the topology
of uniform convergence on compact sets. It is easy
to see that the limit reaches its maximum
on any $B_r$, hence by maximum principle
this limit is constant. On the other hand,
$\sup \gamma^n f$ stays the same, hence
the limit is the same for all subsequences
of $\gamma^n f$. We obtain that $\gamma^n f$
converges to a constant.

The set $Z:= \tilde M_c\backslash\tilde M$
is by construction compact and fixed by
$\Gamma$. Therefore, $\sup_Z |\gamma^n f| = \sup_Z |f|$
and $\inf_Z |\gamma^n f| = \inf_Z |f|$ for any
holomorphic function $f$. Since 
$\gamma^n f$ converges to a constant,
we obtain that $\inf_Z |f|=\sup_Z |f|$.
This implies that all globally defined holomorphic
functions take the same value in all points of $Z$.
The variety $\tilde M_c$ is Stein, and therefore
for any pair $x, y \in \tilde M_c$, $x\neq y$ there
exists a global holomorphic function taking distinct
values at $x, y$. Therefore, $Z$ is just one point.

We obtained that
$\tilde M_c$ is a one-point compactification
of $\tilde M$. This proves \ref{_filling_Theorem_}.
\endproof

\hfill

\proposition\label{_gamma_eigenva_Proposition_}
In assumptions of \ref{_filling_Theorem_},
denote the point $\tilde M_c\backslash\tilde M$ by  $z$.
Let $\goth m$ be the ideal of $z$, and $\gamma$
the monodromy generator, which acts on $\tilde M_c$.
Then  $\gamma$ acts with eigenvalues strictly smaller 
than $1$ on the cotangent space $\goth m/\goth m^2$. 

\hfill

\noindent{\bf Proof:}
It will be enough to carry on the argument on the tangent space. In
the smooth situation, if $\gamma$ fixes a point $x$ and maps an open
neighborhood $U$ of $x$ into a compact neighborhood 
of $x$ contained in $U$, then $\gamma$ will act on $T_xU$ as
desired, as Schwarz lemma implies. 
In general, we already know that $\lim\gamma^n f=const$
for any holomorphic function $f$.
 Now let $v$ be an eigenvector of $d_x\gamma\in \mathrm{End}T_xU$
with eigenvalue $\lambda$, such that $d_xf(v)\neq 0$. Then
$d_x(\gamma^n f)(v)=\lambda^n d_xf(v)$. As  $\gamma^n f$ converges
to a constant, $\lambda^n d_xf(v)$ converges to $0$, 
and hence $\mid \lambda\mid < 1$.
\endproof

\subsection{Holomorphic flow and $\Z$-equivariant maps}

By \ref{_gamma_eigenva_Proposition_}, 
the formal logarithm of $\gamma$
converges, and $\gamma$ is obtained by integrating
a holomorphic flow.

Using an argument similar to \cite{ve96} 
we show that $\gamma$ acts with finite Jordan blocks on
the formal completion $\hat \calo_z$  of the local ring $\calo_z$ 
of germs of analytic functions in $z\in \tilde M_c$. But more is true: 
the formal solution of the ODE is in fact analytic, and hence 
$\gamma$ acts with finite Jordan blocks on $\hat \calo_z$.

\hfill

\theorem
Let $\nu$ be a holomorphic flow on a Stein variety $S$, 
acting with eigenvalues strictly smaller than $1$ on the 
tangent space $T_sS$ for some $s\in S$. Then there exists a sequence $V_n\subset\calo_s S$ of 
finite-dimensional subspaces in the ring of germs of analytic functions  on $S$ such that the $s$-adic completion of $\oplus V_n$ equals the completion 
$\hat \calo_s S$ of $\calo_s S$ and all of the $V_n$ are preserved by $\nu$ (which 
acts by linear transformations of $V_n$).

\hfill

\noindent{\bf Proof:}
Taking an appropriate embedding of $S$, we may assume that $S$ 
is an analytic subvariety of an open ball $B$ in
$\CC^n$. Since $S$ is Stein, the natural restriction map
$T_sB\rightarrow T_sB|_{S}$ is surjective, therefore $\nu$
can be extended to a holomorphic flow on $B$. As a
consequence, it is enough to prove the theorem assuming
$S$ smooth. 

We take $S$ to be an open ball in $\CC^n$ and $s=0$. By the Poincar\'e-Dulac theorem (cf. \cite[\S 23, p. 181]{arnold}), after a biholomorphic change of coordinates, $\nu$ can be given the following normal form:
$$
\nu(x_i)=\la_i x_i +P(x_{i+1},\ldots,x_n)
$$
where $P(x_{i+1},\ldots,x_n)$ are resonant polynomials\footnote{We recall, \emph{loc. cit.}, that 
an $n$-tuple $\la=(\la_1,\ldots,\la_n)$ is called resonant if there exists an integral relation $\la_s=(\la, m):=\sum m_i\la_i$, with $m=(m_1,\ldots,m_n)$, all $m_i\geq 0$ and $\mid m\mid:=\sum m_i\geq 2$. Correspondingly, if the eigenvalue $\la$ of an operator $A$ is resonant and $x^m=x_1^{m_1}\cdots x_n^{m_n}$ is a monomial, $x_i$ being the coordinates of a vector $x$ in a fixed eigenbasis $\{e_i\}$, we say that the vector-valued monomial $x^me_s$ is  resonant if $\la_s=(m,\la)$, $\mid m\mid \geq 2$.} corresponding to the  
eigenvalues $\la_i$.

Let $V_{\la_i}$ be the spaces generated by the coordinate functions $x_i$ corresponding to the eigenvalue $\la_i$ and resonant monomials. By construction, $\nu$ preserves each $V_{\la_i}$. Taking a set $V_{\la_{i_1}}, \ldots , V_{\la_{i_k}}$, we see that $\nu$ also preserves 
the space $V_{\la_{i_1}, \ldots , \la_{i_k}}$ generated by 
$q_1,\ldots, q_k$, with $q_j\in V_j$. But a completion of  all 
these $V_{\la_{i_1}, \ldots , \la_{i_k}}$ is, by construction, equal with $\hat\calo_s S$.

\endproof

\hfill

Picking enough holomorphic functions in these
finite-\-di\-men\-si\-o\-nal ei\-gen\-spaces, we find an embedding
$\tilde M_c\hookrightarrow \CC^n$, such that the 
monodromy $\Gamma$ acts on $\CC^n$ equivariantly. 
Then $M$ is embedded to $(\CC^n\backslash 0) / \Gamma$,
where $\Gamma$ acts linearly with eigenvalues strictly smaller than $1$.

We shall call such a quotient a \emph{Hopf manifold};
to avoid confusion with class 0 Hopf surfaces of Kodaira -
"linear Hopf manifold". It can be viewed as a generalization in 
any dimension of the class 0  Hopf surface. Clearly, such
a Hopf  manifold is LCK (\cite{go}) but, in general, not Vaisman. 

Using Strenberg's normal form of holomorphic contractions 
(cf. \cite{sternberg}), Belgun also showed that Hopf manifolds are LCK and 
classified the Vaisman ones (they are precisely those for which the normal form of the contraction is resonance free), relating them to Sasakian geometry 
(see \cite{belgun_video}).

As a result, we obtain the following:

\hfill

\theorem\label{_koda_embe_Theorem_}
Any  compact LCK manifold $M$ with potential, $\dim M \geq 3$, 
can be embedded to a linear Hopf manifold.

\endproof

\remark
The converse is also true, of course.

\hfill

This gives a new version of the Kodaira-type theorem 
on immersion of Vaisman manifolds to Vaisman-type 
Hopf manifolds (cf. \cite{ov2}). 
In fact, a much stronger result is obtained. 

\hfill

\theorem\label{vai}
Let $M$ be a Vaisman manifold. Then $M$ can be embedded to a 
Vaisman-type Hopf manifold $H = (\C^n\backslash 0) /\langle A\rangle$,
where $A$ is a diagonal linear operator on $\C^n$
with all eigenvalues $<1$.

\hfill

\noindent{\bf Proof:}
Let $S_a:= \phi^{-1}(a)$  be the level set of 
$\phi$, and $V$ the space of CR-holomorphic functions
on $S$. Using the solution of $\bar\6$-Neumann problem
(\cite{ma}), we identify the $L^2$-completion of $V$
with the space of $L^2$-integrable holomorphic functions
on the pseudoconvex domain bounded by $S_a$.

Let $X:=\theta^\sharp$ be the 
holomorphic vector field dual to the Lee form.
According to the general results from
Vaisman geometry (\cite[p. 37]{drag}), $X$ is
holomorphic. Therefore, $X$ acts on $V$
in a natural way. Denote by $g$ the natural
metric on $S_a$. Since $X$
acts on $\tilde M$ by homotheties
(\cite[p. 37]{drag}), we have $X(g) = cg$. 
Therefore, $X$ acts on $V$ as a self-dual operator.
Then, $X$ is diagonal on any finite-dimensional
subspace in $V$ preserved by $X$. This is true
for the spaces $V_{\la_{i_1}, \ldots , \la_{i_k}}$
used to construct the embedding $\tilde M_c \arrow \C^n$.

\hfill\endproof

\hfill

After this work was completed, professor Hajime Tsuji
communicated to us the reference to 
\cite{_Kato2_}. In this work, Masahide Kato obtains 
a sufficient condition for a compact complex manifold 
$Z$ of $\text{dimension}\geq 4$ to dominate
bimeromorphically a subvariety of a Hopf manifold,
in terms of a certain effective divisor on $D$
and a flat line bundle on $Z$.

Also, in \cite{_Kato1_}, the following theorem was obtained. Let
$X$ be a complex space with a fixed point $z$, equipped
with a holomorphic contraction, that is, an invertible 
morphism $\psi:\; X \arrow X$ such that for any sufficient
small neighbourhood $U$ of $z$, $\psi^k(U)$ lies
inside $U$ for $k$ sufficiently big, and for all $x\in X$,
\[ \lim\limits_{k \arrow \infty} \psi^k(x)=z.
\]
Then $(X\backslash z)/\langle \psi\rangle$
can be embedded into a Hopf manifold.

\hfill

{\bf Acknowledgements:} Misha Verbitsky is grateful
to  D. Kaledin, D. Kazhdan, D. Novikov and P. Pushkar for interesting
discussions, and D. Novikov for a lecture on normal forms and 
Poincar\'e-Dulac theorem. Also, a gratitude to IMAR in Bucharest for 
its hospitality during one week in May 2004, when this research 
was initiated. Both authors thank the referee for very useful comments which improved the presentation.

{\small

\noindent {\sc Liviu Ornea\\
University of Bucharest, Faculty of Mathematics, \\14
Academiei str., 70109 Bucharest, Romania.}\\
\tt Liviu.Ornea@imar.ro, \ \ lornea@gta.math.unibuc.ro

\hfill

\noindent {\sc Misha Verbitsky\\
University of Glasgow, Department of Mathematics, \\15
  University Gardens, Glasgow, Scotland.}\\
{\sc  Institute of Theoretical and
Experimental Physics \\
B. Cheremushkinskaya, 25, Moscow, 117259, Russia }\\
\tt verbit@maths.gla.ac.uk, \ \  verbit@mccme.ru 
}


{\scriptsize
\begin{thebibliography}{100}

\bibitem[AS]{andreotti_siu} 
A. Andreotti and Y.T. Siu, \emph{Projective embeddings of
pseudoconcave spaces}, Ann. Sc. Norm. Sup. Pisa {\bf 24} (1970),
231--278.

\bibitem[Ar]{arnold} V.I. Arnold, Geometrical methods in the theory of 
ordinary differential equations, Grundlehren der Mathematischen Wissenschaften 
250, Springer Verlag, 1983.



\bibitem[B1]{belgun} F.A. Belgun, {\em On the metric structure of
non-K{\"a}hler complex surfaces}, Math. Ann. {\bf 317} (2000),
1--40.

\bibitem[B2]{belgun_video} F.A. Belgun, \emph{Hopf manifolds and
Sasakian structures}, Lecture given at  CIRM 2002, ``G\'eom\'etrie
des vari\'et\'es de petites dimensions et g\'eom\'etries
sp\'eciales'', video available at
www.cirm.univ-mrs.fr/videos/index.php

\bibitem[Bl]{bl} D.E. Blair, Riemannian geometry of contact and symplectic manifolds, Progress in Math. {\bf 203}, Birkh\"auser, Boston, Basel 2002.

\bibitem[Dr]{sorin} S. Dragomir, \emph{On pseudo-Hermitian immersions between strictly pseudoconvex CR manifolds}, Amer. J. Math. 117 (1995), no. 1, 169--202.

\bibitem[DO]{drag} 
S.  Dragomir and L.  Ornea,  Locally conformal
K{\"a}hler
geometry, Progress in Math. {\bf 155},   Birkh{\"a}user, Boston, Basel, 1998.

\bibitem[GO]{go} 
P. Gauduchon and L. Ornea,
{\em Locally conformally K{\"a}hler metrics on Hopf surfaces},
Ann. Inst. Fourier  {\bf 48} (1998), 1107--1127.

\bibitem[KO]{_Kamishima_Ornea_}
Y. Kamishima, L. Ornea,
{\em Geometric flow on compact locally conformally K\"ahler manifolds},  Tohoku Math. J. {\bf 57}  (2005), 201--221.

\bibitem[Ka1]{_Kato1_}
Kato, Ma. {\em Some remarks on subvarieties of Hopf manifolds}, 
A Symposium on Complex Manifolds (Kyoto, 1974). 
S\^urikaisekikenky\^usho K\'oky\^uroku No. 240 (1975), 64--87.

\bibitem[Ka2]{_Kato2_} 
Kato, Ma. 
{\em On a characterization of submanifolds of Hopf manifolds}, 
Complex analysis and algebraic geometry, pp. 191--206. 
Iwanami Shoten, Tokyo, 1977.

\bibitem[MD]{ma} 
G. Marinescu and  T.-C. Dinh, \emph{On the
compactification of hyperconcave ends and the theorems of Siu-Yau and Nadel} math.CV/0210485, v2.

\bibitem[MY]{my} G. Marinescu and N. Yeganefar, \emph{Embeddability of some strongly pseudoconvex CR manifolds}, math.CV/0403044.


\bibitem[OV1]{ov1} L. Ornea and M. Verbitsky, \emph{Structure theorem for compact Vaisman manifolds}, Math. Res. Lett. {\bf 10} (2003), 799--805.

\bibitem[OV2]{ov2} L. Ornea and M. Verbitsky, \emph{An immersion theorem for
compact Vaisman manifolds},  Math. Ann. {\bf 332}  (2005),  no. 1, 121--143. 

\bibitem[Re]{renaud} J. Renaud, \emph{Classes de vari\'et\'es
localemant conform\'ement K\"ahleriennes non K\"ahleriennes},
C.R. Acad. Sci. Paris {\bf 338} (2004), 925--928.

\bibitem[Ro]{rossi}
H. Rossi, \emph{Attaching analytic spaces to a analytic space along
a pseudo-convex boundary}, Proc. Conf. Complex Manifolds
(Minneapolis), Springer Verlag, 1965, 242--256.

\bibitem[S]{sternberg} S. Sternberg, \emph{Local contractions and a theorem of Poincar\'e}, Amer. J. Math. {\bf 79} (1957), 809--824.

\bibitem[Va]{_Vaisman:Dedicata_}
I. Vaisman,  {\em  Generalized Hopf
manifolds}, Geom. Dedicata {\bf 13} (1982), no. 3, 231--255.

\bibitem[V1]{ve96} M. Verbitsky, \emph{Desingularization of singular
hyperk\"ahler varieties}, II, alg-geom/9612013.

\bibitem[V2]{ve} M. Verbitsky, \emph{Theorems on the vanishing of cohomology for locally conformally hyper-Kähler manifolds} (Russian).  Tr. Mat. Inst. Steklova  {\bf 246}  (2004),  Algebr. Geom. Metody, Svyazi i Prilozh., 64--91;  translation in  Proc. Steklov Inst. Math.  2004,  {\bf 246}, 54--78. Also available as math.DG/0302219.
\end{thebibliography}
}
\end{document}